# A Novel Discrete Adjoint-based Level Set Topology Optimization Method in B-spline Space

Hao Deng*, Kazu Saitou

Department of Mechanical Engineering, University of Michigan, Ann Arbor, MI 48109

*Corresponding author. Email: denghao@umich.edu

**Abstract**

This paper presents a novel computational scheme for sensitivity analysis of the velocity field in the level set method using the discrete adjoint method. The velocity field is represented in B-spline space, and the adjoint equations are constructed based on the discretized governing equations. The key contribution of this work is the demonstration that the velocity field in the level set method can be entirely obtained from the discrete adjoint method. This eliminates the need for shape sensitivity analysis, which is commonly used in standard level set methods. The results demonstrate the effectiveness of the approach in producing optimized results for stress and linearized buckling problems. Overall, the proposed method has the potential to simplify the way in which topology optimization problems using level set methods are solved, and has significant implications for the design of a broad range of engineering applications.

## 1. Introduction

Continuum topology optimization has emerged as a crucial tool for determining optimal shapes that maximize performance, subject to specified design constraints. Since the initial development of the Homogenization Design Method (HMD) [1] three decades ago, numerous topology optimization methods have been proposed. These include the density-based method [2], level set method [3], BESO method [4], and more recently, approaches such as moving morphable components (MMC) [5-7], moving morphable voids (MMV) [8], and the geometry projection method [9-11].

The Level set method (LSM) is a class of shape and topology optimization method, where the topology of a structure is described implicitly using the level set function $\phi$, with the boundary defined by $\{\phi = 0\}$. A comprehensive literature review for level set methods in topology optimization can be found in Ref [12]. In classical LSM, the boundary is updated through evolving the implicit function $\phi$ based on Hamilton-Jacobi equation in which the optimization problem needs to be transformed in an unconstrained problem by Lagrangian formulations. During the optimization, the Lagrangian multipliers are gradually updated based on a certain strategy to improve convergence. The evolution of the implicit level set function is driven by normal velocity of the boundary $\{\phi = 0\}$ based on natural velocity extension [13] or the fast-marching method [14]. Besides the classical LSM formulation, the parametric level-set method has drawn great attention in recent years. Wang et al [15] proposed an Radial Basis Function (RBF)level set optimization method to transform the Hamilton-Jacobi equation into a system of ordinary differential equations (ODEs) based on a collocation formulation. Since then, several parametric level-set methods [16-19] with different basis functions for various physical problems have been proposed. In a velocity field level-set (VFLS) method [20, 21] recently proposed by Wang et al, the parametric boundary normal velocity field is defined to derive the level set function, instead of the parametric level set function deriving the boundary normal velocity as conventionally formulated. The velocity field is controlled by prescribed basis functions and velocity design variables at given points in the design domain. VFLS method provides an effective way to handling multiple constraints and enables the use of the standard mathematical programming algorithms. In [22], Wang et al incorporated the topological derivative concept into VFLS method to enable the automatic nucleation of interior holes.

For classical LSM for topology optimization, the continuous adjoint method is widely used to compute sensitivities [13]. One of the major obstacles of continuous adjoint method is the sensitivity expressions are discontinuous at the nodes and edges of finite elements. In general, interpolation and smoothing techniques are needed to avoid discontinuity with a price of reduced the accuracy of the resulting sensitivity. Such interpolation and smoothing may even lead to divergence for high move limits as demonstrated by Kambampati et al [23]. Compared to the continuous adjoint method that follows the "differentiate-then-discretize" scheme [24], the discrete adjoint method follows the "discretize-then-differentiate" scheme, where the partial differential equations are first discretized using finite element method. The discretized equations are then differentiated to obtain adjoint equations based on augmented functional equations. While the discretize adjoint method has been widely used in density-based topology optimization methods [25], Kambampati et al [23] introduced the discrete adjoint method into the classical LSM for the first time using a semi-analytical formulation, where the boundary is perturbed implicitly to obtain level set sensitivity (*i.e.*, velocity) based on finite difference approximation. Reference [26] presents a comprehensive literature review of existing level set methods for topology optimization, providing a detailed comparison and analysis of various approaches. Dunning et al [27] introduced linear buckling constraints for level-set topology optimization. In their approach, the velocity function is defined as a weighted sum of the shape sensitivities for the objective and constraint functions. R. Picelli et al [28] proposed a level set method to solve minimum stress and stress-constrained shape and topology optimization problems, where the shape sensitivity function is derived and a computational procedure based on a least squares interpolation approach is devised in order to compute sensitivities at the boundaries. Utilizing the level set method to address topology optimization problems involving both stress and buckling constraints is challenging and, to the best of the author's knowledge, has not been reported previously. Nevertheless, considering both constraints is crucial for real-world engineering applications.

This paper presents an analytical formulation for the discrete adjoint sensitivity in level set method (LSM)-based topology optimization, establishing a pioneering connection between the discrete adjoint method in classical LSM and the normal boundary velocity. The proposed analytical formulation enables the derivation of the normal boundary velocity of a level set from the discrete adjoint method analytically, thereby bypassing shape sensitivity analysis. In this paper, the velocity field is described in B-spline space, with similar parameterization methods for topology optimization found in Refs. [29-31]. The effectiveness of the proposed formulation is demonstrated through VFLS-based topology optimization using B-spline basis, considering stress and local buckling constraints.

The remainder of the paper is organized as follows: Section 2 presents the mathematical formulation of B-spline VFLS based on the discrete adjoint method. Section 3 showcases two typical topology optimization problems incorporating local stress and buckling constraints. Section 4 concludes the paper, discussing possible future directions.

## 2. B-spline Velocity Field Level Set Topology Optimization Method

## 2.1 Definition of B-spline velocity field level set (VFLS) function

In LSM, the material domain $\Omega$ at time $t$ during the optimization within fixed design domain $D$ can be represented by a signed distance function $\Phi(x, t)$:

$$\begin{cases} \Phi(x,t) > 0 \text{ if } x \in \Omega \setminus \partial\Omega \\ \Phi(x,t) = 0 \text{ if } x \in \partial\Omega \cap D \\ \Phi(x,t) < 0 \text{ if } x \in D \setminus \Omega \end{cases} \quad (1)$$

The shape and topology of material domain $\Omega$ evolve over a course of time through moving its boundary $\partial\Omega$, which is realized by solving the Hamilton-Jacobi equation with velocity $V^n$ in the normal direction of $\partial\Omega$:

$$\frac{\partial \Phi}{\partial t} - V^n |\nabla \Phi| = 0 \quad (2)$$

For $D \subset \mathbb{R}^2$, we can construct normal velocity field $V^n(x, t)$ using a piece-wise function as follows,

$$V^n(x, y, t) = \sum_{j=1}^{n_K} V_j(t) p_j(x, y) \quad (3)$$

where $V_j$ is design variables. $n_K$ is the number of design variable points (velocity knots). In this paper, the velocity knots are selected as the centroid of FEM element. A special piecewise function is utilized here to work as basis function,

$$p_j(x) = \begin{cases} 1 & \|x - x_j\|_\infty \leq r \\ 0 & \|x - x_j\|_\infty > r \end{cases} \quad (4)$$

where $x_j$ denotes the centroid of the $j$ element. For 2D square mesh, $r = 0.5L$ is selected for computation, where $L$ is element length. Based on this special piecewise basis, the velocity field inside each element is a constant and only depends on the velocity at the centroid of element. This special piecewise function brings a simplification for sensitivity computation. According to the velocity field level set method, the movement of boundary depends on the boundary velocity. Due to the piecewise function introduced in Eq. (4), the boundary velocity inside one element is a constant and only depends on the velocity value at the centroid of current element (velocity knot), which is independent from neighborhood velocity knots. Thus, the change of element density (or volume) is only dependent on the current element velocity knot. This property can be described in a mathematical way as shown in Eq. (13) in section 2.2. The value at velocity knots can be described in B-spline space as follows,

$$V(x, y, t) = \sum_{k=0}^{n_x} \sum_{l=0}^{n_y} B_{k,p}(x) B_{l,q}(y) b_{k,l}(t) \quad (5)$$

where $B_{k,p}: \mathbb{R} \to \mathbb{R}$ $B_{l,q}: \mathbb{R} \to \mathbb{R}$ ($k = 0,1, \ldots, n_x; l = 0,1, \ldots n_y$) are B-spline basis functions defined by knot vectors in $x$ and $y$ directions, respectively, $p$ and $q$ are the degrees of the B-spline basis functions, and $b_{k,l}(t)$ is a B-spline coefficient [32]. Once the velocity field is obtained by Eq. (3) and Eq. (5) for a given value of $b_{k,l}(t)$ at time $t$ during optimization, level set function $\Phi(x, t)$ can be updated through Eq. (2) using, for example, the upwind difference scheme [14] and the method of moving asymptotes (MMA) algorithm [33]. Then, the re-initialization step is implemented to avoid numerical deterioration of the level set function. More implementation details can be found in Ref [20].

## 2.2 Sensitivity analysis based on discrete adjoint method

The analytical derivation of the normal velocity field $V^n$ from the discrete adjoint method is presented in two stages. Initially, we discuss the derivation of the discrete adjoint sensitivity for the density-based method, as outlined in reference [34]. Subsequently, we employ the chain rule to establish a connection between the sensitivity concerning B-spline velocity coefficients and the discrete adjoint sensitivity with respect to density, which was determined in the first stage

**Step 1:** The discretized governing equation for linear elastic problem can be written as follows,

$$\boldsymbol{\Psi} = \boldsymbol{K}\boldsymbol{u} - \boldsymbol{f} = \boldsymbol{0} \tag{6}$$

where $\boldsymbol{K}$ is the stiffness matrix, $\boldsymbol{u}$ is the displacement vector, and $\boldsymbol{f}$ is the nodal force vector. Note that any other constraint equations for physical problems are written as follows,

$$\boldsymbol{\mathcal{H}} = \boldsymbol{0} \tag{7}$$

For a given objective function $F$, an augmented Lagrangian function is defined as follows,

$$G = F + \boldsymbol{\psi}^T \boldsymbol{\Psi} + \boldsymbol{\kappa}^T \boldsymbol{\mathcal{H}} \tag{8}$$

where $\boldsymbol{\psi}$ and $\boldsymbol{\kappa}$ are Lagrange multipliers. The sensitivity of augmented Lagrangian function $G$ with respect to discretized density field $\boldsymbol{\rho}$ within D can be expressed as follows,

$$\frac{\partial G}{\partial \boldsymbol{\rho}} = \frac{\partial F}{\partial \boldsymbol{\rho}} + \boldsymbol{\psi}^T \frac{\partial \boldsymbol{\Psi}}{\partial \boldsymbol{\rho}} + \boldsymbol{\kappa}^T \frac{\partial \boldsymbol{\mathcal{H}}}{\partial \boldsymbol{\rho}} \tag{9}$$

where Lagrange multiplier should satisfy,

$$\boldsymbol{\psi}^T \frac{\partial \boldsymbol{\Psi}}{\partial \boldsymbol{\rho}} + \boldsymbol{\kappa}^T \frac{\partial \boldsymbol{\mathcal{H}}}{\partial \boldsymbol{\rho}} + \frac{\partial F}{\partial \boldsymbol{\rho}} = \boldsymbol{0} \tag{10}$$

The above derivation is usually referred to as the discrete adjoint method.

**Step 2:** Based on the chain rule,

$$\frac{\partial F}{\partial b_{k,l}} = \sum_{i=1}^{n_e} \frac{\partial F}{\partial \rho_i} \frac{\partial \rho_i}{\partial b_{k,l}} \tag{11}$$

The first term $\frac{\partial F}{\partial \rho_i}$ in the right-hand side of Eq. (11) can be obtained by Eq. (9). Based on the chain rule, the second term $\frac{\partial \rho_i}{\partial b_{k,l}}$ can be written as,

$$\frac{\partial \rho_i}{\partial b_{k,l}} = \sum_{j=1}^{n_v} \frac{\partial \rho_i}{\partial V_j^n} \frac{\partial V_j^n}{\partial b_{k,l}} \tag{12}$$

Here, $V_j^n$ and $\rho_i$ denote the normal velocity at velocity point $j$ in $\partial \Omega$ and the density at element $i$ in D, respectively, with $n_e$ representing the number of finite elements in D and $n_v$ signifying the number of velocity points in $\partial \Omega$. The first term $\frac{\partial \rho_i}{\partial V_j^n}$ on the right-hand side of Eq. (12) indicates the change in the volume fraction (density) of element $i$ resulting from the alteration in the normal velocity point $j$. Unlike Kambampati et al. [23], where implicit perturbation near the boundary point is employed to numerically approximate the second term $\frac{\partial \rho_i}{\partial V_j^n}$ using the finite difference method, in this paper we analytically derive $\frac{\partial \rho_i}{\partial V_j^n}$ without numerical approximation. We can express this as:

$$\frac{\partial \rho_i}{\partial V_j^n} = \delta_{ij} \frac{\partial \rho_i}{\partial V_j^n} \tag{13}$$

where $\delta_{ij} = 1\ (if\ i = j), and\ \delta_{ij} = 0\ (if\ i \neq j)$. Note that the above equation is based on the piecewise property, as demonstrated in Eq. (3) and (4), which means the change in each element density (or volume) is solely dependent on the current element velocity. The above equation can be reformulated as:

$$\frac{\partial \rho_i}{\partial V_j^n} = \begin{cases} \frac{\partial \rho_i}{\partial V_i^n}\ (j = i) \\ 0\ (otherwise) \end{cases} \tag{14}$$

Since $\rho_i = H(\Phi_i)$ where $H(\cdot)$ is Heaviside function [3], the chain rule yields:

$$\frac{\partial \rho_i}{\partial V_j^n} = \delta_{ij} \frac{\partial \rho_i}{\partial V_i^n} = \delta_{ij} \frac{\partial H(\Phi_i)}{\partial V_i^n} = \delta_{ij} \delta(\Phi_i) \cdot \frac{\partial \Phi_i}{\partial V_i^n} \tag{15}$$

where $\delta(x) \equiv \frac{\partial H(x)}{\partial x}$ is Dirac function [3]. In order to obtain $\frac{\partial \Phi_i}{\partial V_i^n}$, Eq. (2) is rewritten as:

$$d\Phi_i - V_i^n |\nabla \Phi_i| dt = 0 \qquad (i = 1,2,\cdots, n_v) \tag{16}$$

where $dt$ is a pseudo time, $\Phi_i \equiv \Phi(x_i, t)$, $V_i^n \equiv V^n(x_i, t)$. Now let $\delta V_i^n$ be a small perturbation of velocity $V_i^n$ and $\delta d\Phi_i$ be the corresponding variation of $d\Phi_i$. Ignoring the effect of small perturbation $\delta d\Phi_i$ on $|\nabla \Phi_i|$, Eq (16) becomes:

$$(d\Phi_i + \delta d\Phi_i) - (V_i^n + \delta V_i^n)|\nabla \Phi_i| dt = 0 \tag{17}$$

which can be simplified as:

$$\delta d\Phi_i - \delta V_i^n |\nabla \Phi_i| dt = 0 \tag{18}$$

Therefore,

$$\frac{\partial d\Phi_i}{\partial V_i^n} = |\nabla \Phi_i| dt \tag{19}$$

which can be rewritten as,

$$\frac{\partial (d\Phi_i/dt)}{\partial V_i^n} = |\nabla \Phi_i| \tag{20}$$

The level set function is usually defined as a signed distance function. To preserve the signed distance property, the following re-initialization step is usually required in each iteration:

$$\frac{\partial \Phi}{\partial t} + \text{sign}(\Phi_0)(|\nabla \Phi| - 1) = 0 \tag{21}$$

where $\Phi_0$ is the level set function before re-initialization. After re-initialization, we have

$$|\nabla \Phi| \approx 1 \tag{22}$$

Based on Eq. (20), we have,

$$\frac{\partial (d\Phi_i/dt)}{\partial V_i^n} = \frac{d(\partial \Phi_i/\partial V_i^n)}{\partial t} = |\nabla \Phi| \approx 1 \tag{23}$$

Therefore,

$$\frac{\partial \Phi_i}{\partial V_i^n} = \int |\nabla \Phi| dt \approx \int 1 \, dt \tag{24}$$

It is important to note that as long as the re-initialization of the level set function during each step can achieve high accuracy, ensuring that $|\nabla \Phi| = 1$, the accuracy of the velocity field sensitivity can be guaranteed. Based on Eq. (23), we deduce:

$$\frac{\partial \Phi_i}{\partial V_i^n} \approx t \tag{25}$$

It is important to note that the pseudo time $t$ does not hold any physical meaning here and can be regarded as 1. This is because $t$ merely serves as a scaling factor for sensitivity $\frac{\partial F}{\partial V_j^n}$, which will not impact the optimization progress if the optimization algorithm can adjust move steps at each iteration, such as the method of moving asymptotes (MMA) [30]. Therefore, by substituting Eq. (25) and Eq. (15) into Eq. (12) and utilizing Eq. (5), we obtain:

$$\frac{\partial \rho_i}{\partial b_{k,l}} = \delta(\Phi_i) \cdot \sum_{j=1}^{n_v} \left( \delta_{ij} B_{k,p}(x_j) B_{l,q}(y_j) \right) \tag{26}$$

Eq. (26) can be further simplified as:

$$\frac{\partial \rho_i}{\partial b_{k,l}} = \delta(\Phi_i) B_{k,p}(x_i) B_{l,q}(y_i) \tag{27}$$

where $(x_i, y_i)$ is the coordinates of the centroid of element $i$. Substituting this to Eq. (11), the sensitivity of function $F$ with respect to B-spline coefficient $b_{k,l}$ can be given as:

$$\frac{\partial F}{\partial b_{k,l}} = \sum_{i=1}^{n_e} \frac{\partial F}{\partial \rho_i} \delta(\Phi_i) B_{k,p}(x_i) B_{l,q}(y_i) \tag{28}$$

It is worth noting that the level set B-spline representation can be regarded as a filter that influences the final features of the optimized shape, similar to the effect of the B-spline filter in density-based methods. B-spline functions provide a smooth and continuous representation of the velocity field, which can effectively control the shape complexity and ensure the generation of high-quality designs.

## 2.3 Topology optimization problem formulation

In the examples in the next section, the B-spline velocity field level set topology optimization is formulated as the volume minimization subject to stress and local buckling constraints. The optimization problem can be described as follows:

$$\begin{aligned} & \min V \\ & s.t. \ \sigma^{PM} \leq \sigma^* \\ & s.t. \ KS[\mu_i]_{(i \in \mathbb{Z})} \leq \mu^* \end{aligned} \tag{29}$$

where $\sigma^{PM}$ and $\sigma^*$ respectively denote the p-norm stress and its upper limit for constraining the maximum stress, and $KS[\mu_i]_{(i \in \mathbb{Z})}$ and $\mu^*$ respectively denote the KS aggregation and its upper bound for guarding against buckling. The reader should refer to [29] and [30] for the mathematical formulation and derivation of the sensitivities for these constraints in detail. While these formulations are given as a function of element density $\rho_i$, the sensitivities with respect to B-spline coefficients $b_{k,l}$ of velocity field can be readily derived based on the mathematical formulation in Eq. (22).

### 2.3.1. Sensitivity of stress constraint based on discrete adjoint method

The discretized linear elastic equations can be written as,

$$\boldsymbol{K}\boldsymbol{u} = \boldsymbol{f} \tag{30}$$

where $\boldsymbol{K}$ is the global stiffness matrix, $\boldsymbol{u}$ is displacement vector, and $\boldsymbol{f}$ is force vector. The global stiffness matrix can be assembled by elemental stiffness matrix $\boldsymbol{K}_{e,i}$ as follows,

$$\boldsymbol{K} = \sum_{i=1}^{n_e} \boldsymbol{K}_{e,i} \tag{31}$$

where $\Sigma$ here represents the element matrix assembly operator. For density-based topology optimization, the element stiffness matrices are modeled as Ersatz material:

$$\boldsymbol{K}_{e,i} = \left(E_{min} + \rho_i(E - E_{min})\right)\boldsymbol{K_0} \tag{32}$$

where $\rho_i$ is the density (volume fraction) of the $ith$ element, $E_{min}$ is the elasticity modulus of a void element, which is a small value to avoid numerical issue, and $\boldsymbol{K_0}$ is the element stiffness matrix for unit elasticity modulus. The maximum von-mises stress can be approximated by p-norm stress [35], which can be written as,

$$\sigma^{PM} = \left(\sum_i^{n_e} \sigma_{vm,i}^p\right)^{\frac{1}{p}} \tag{33}$$

where $\sigma_{vm,i}$ is the von Mises stress of the $ith$ element. As described in Section 2.2, the augmented Lagrangian function $\mathcal{L}_p$ can be formulated as follows,

$$\mathcal{L}_p = \sigma^{PM} + \boldsymbol{\lambda}_p^T(\boldsymbol{K}\boldsymbol{u} - \boldsymbol{f}) \tag{34}$$

The adjoint vector $\boldsymbol{\lambda}_p$ can be obtained through solving $\frac{\partial \mathcal{L}_p}{\partial u} = 0$. Therefore, the sensitivity of p-norm stress with respect to elemental density is given by,

$$\frac{\partial \mathcal{L}_p}{\partial \rho_i} = \frac{\partial \sigma^{PM}}{\partial \rho_i} + \boldsymbol{\lambda}_p^T \frac{\partial \boldsymbol{K}}{\partial \rho_i}\boldsymbol{u} \tag{35}$$

The detailed derivation of p-norm stress sensitivity can be found in Ref [28].

### 2.3.2 Sensitivity of buckling constraint based on discrete adjoint method

In this part, the sensitivity of fundamental buckling load factor (BLF) with respect to element density is derived. Similarly, the elemental stiffness matrix can be represented by Ersatz material model as shown in Eq. (32). The fundamental buckling load factor $\lambda_f$ can be formulated as Rayleigh quotient [36],

$$\lambda_f(\boldsymbol{\rho},\boldsymbol{u}) = \min\left(-\frac{v^T K[\rho]v}{v^T G[\rho,u]v}\right); \ (v \in \mathbb{R}^n, v \neq \boldsymbol{0}) \tag{36}$$

where $\boldsymbol{G}$ is the global stress stiffness matrix, $\boldsymbol{u}$ is the displacement vector. $\boldsymbol{K}$ is the linear stiffness matrix. The general procedure of linearized buckling analysis is as follows, a) Define a reference load vector $\boldsymbol{f}$ b) Compute the equilibrium displacement $\boldsymbol{u} = \boldsymbol{K}[\rho]^{-1}\boldsymbol{f}$ c) Set up the stress stiffness matrix $\boldsymbol{G}[x,\boldsymbol{u}(x)]$ d) compute the buckling load and corresponding buckling mode $(\lambda_i, \boldsymbol{\varphi}_i)$ using the following equations,

$$(K(\rho) + \lambda G[\rho, u(\rho)])\varphi = 0, \quad \varphi \neq 0 \tag{37}$$

The buckling modes are normalized such that $\varphi_i^T K[\rho]\varphi_j = \delta_{ij}$. Note that the sensitivity of the $ith$ eigenvalue $\lambda_i$ with respect to element density $\rho_e$ is expressed as [37],

$$\frac{\partial \lambda}{\partial \rho_e} = \varphi_i^T \left(\frac{\partial K}{\partial \rho_e} + \lambda_i \frac{\partial G}{\partial \rho_e}\right)\varphi_i - \lambda_i z_i^T \frac{\partial K}{\partial \rho_e} u \tag{38}$$

where $z_i$ can be obtained through solving adjoint equations,

$$Kz_i = \varphi_i^T(\nabla_u G)\varphi_i \tag{39}$$

For local buckling constraint, an aggregation function is applied to generate a single constraint. Here, the KS function [38] is implemented as follow,

$$KS[\mu_i] = \mu_1 + \frac{1}{\gamma} ln\left(\sum_{i=1}^{q} e^{\gamma(\mu_i - \mu_1)}\right) \quad (i \in \mathbb{Z}) \tag{40}$$

Note that $\mu_i = \frac{1}{\lambda_i}$, where the degree of smoothness is governed by the aggregation parameter $\gamma$. $\mathbb{Z}$ is the set of interested $\mu_i$. The value obtained by KS function produce an upper bound of $\max_{i \in \mathbb{Z}} |\mu_i|$. The range of parameter $\gamma$ should be $[1, \infty]$. The first order derivative of KS function $KS[\mu_i]$ with respect to density $\rho_e$ can be written as,

$$\frac{\partial KS[\mu_i]}{\partial \rho_e} = \frac{\sum_{i=1}^{q}\left(e^{\gamma(\mu_i - \mu_1)}\frac{\partial \mu_i}{\partial \rho_e}\right)}{\sum_{i=1}^{q} e^{\gamma(\mu_i - \mu_1)}} \tag{41}$$

The detailed sensitivity derivation and implementation can be found in Ref. [36]. Note that the parameters $p$ and $\gamma$ are set as 8 and 50 in this paper, respectively.

## 4. Numerical Examples

### 4.1 Compressed square design

The first example examines a compressed square design, as depicted in Fig. 1. The square is discretized into $300 \times 300$ finite elements, each with a unit element length. The material's elastic modulus is $E = 1$, and Poisson's ratio is set to $\nu = 0.3$. The bounds for the velocity design variables are chosen as $\pm 0.2$, while the moving limit of the MMA algorithm is 0.2. B-spline knots are uniformly distributed across the design domain with a fixed interval $\Delta$, as shown in Fig. 1. The initialization of the level set function can be seen in Fig. 2. The optimization process terminates when the relative difference in the target function between two consecutive iterations falls below $10^{-4}$. The upper and lower bounds of the normal velocity are set to be within the range [-0.1, 0.1]. For the 2D B-spline basis function, the parameters p and q are selected as p=3 and q=2, respectively.

#### 4.1.1 Volume minimization with p-norm stress constraint

The objective of this example is to minimize the volume fraction while considering only the p-norm stress constraint (i.e., no buckling constraint). The load is uniformly distributed across the middle of the top surface with a total magnitude of $|F| = 1$. The upper bound of p-norm stress, $\sigma^*$, *is set at* $\sigma^*$=1.3. In this simple example, the B-spline knot span length, $\Delta$, is chosen as $\Delta = 0.01L$, and the resulting optimized structure is shown in Fig. 3. The lowest volume fraction ($V = 0.15$) is achieved when $\Delta = 0.01L$, where L represents the length of the square's sides. The initial von Mises stress distribution for the initial design is illustrated in Fig. 2(b), while the final stress distribution for the optimized structure is depicted in Fig. 3(b). The maximum stress decreases from 1.1 to 0.8, and the volume fraction drops from approximately 0.6 at the beginning to 0.15 after 500 iterations. Initially, the von Mises stress distribution for the design is presented in Fig. 2(b). As the optimization process progresses, the stress distribution evolves, ultimately resulting in a more uniform distribution within the optimized structure, as shown in Fig. 3(b). The maximum stress experienced by the structure is reduced from 1.1 to 0.8, indicating that the optimization has led to a more effective design. Furthermore, the volume fraction decreases significantly over the course of 500 iterations, going from an initial value of approximately 0.6 to the final value of 0.15. This reduction in volume fraction highlights the optimization's success in minimizing material usage while simultaneously ensuring that stress constraints are satisfied.

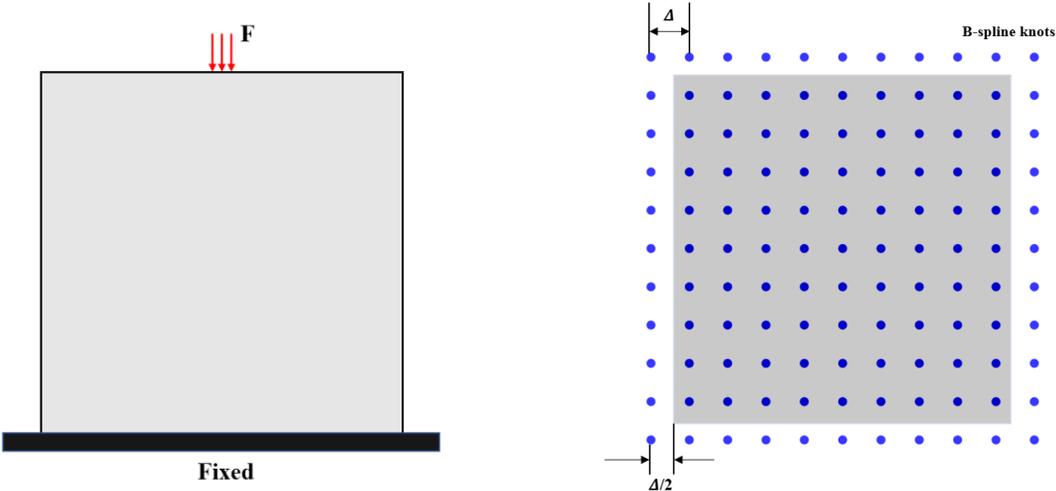

Figure 1. Compressed square design.

(a) material layout　　　　　　　　　　　　　　(b) von Mises stress distribution

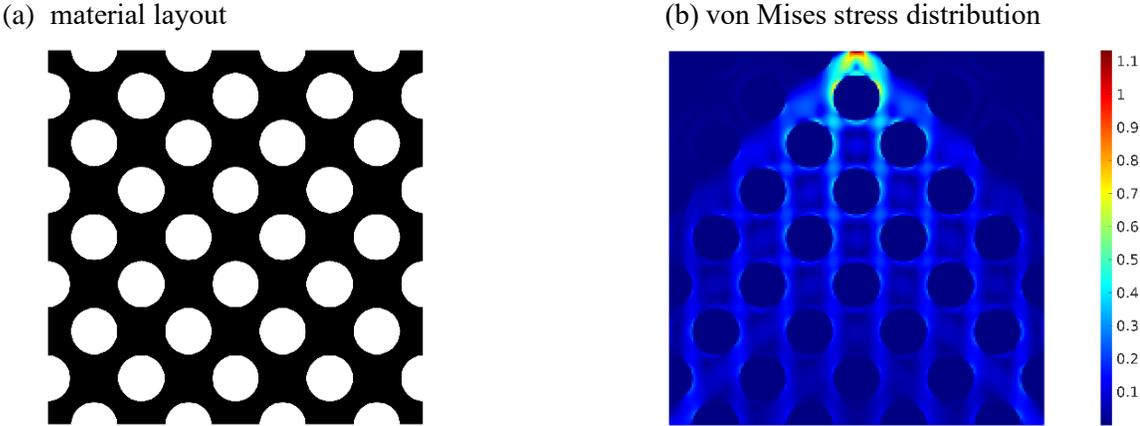

Figure 2. Level set initialization.

(a) material layout　　　　　　　　　　　　　　(b) von Mises stress distribution

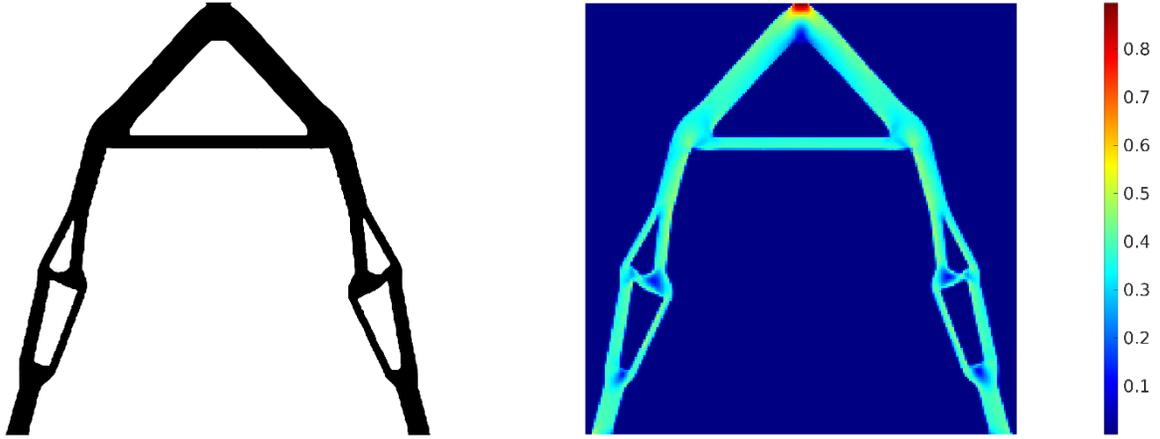

Figure 3. Optimized result: (a) material layout (b) von Mises stress distribution ($V = 0.15, \Delta = 0.01L$).

### 4.1.2 Volume minimization with buckling constraint

The objective of this example is to minimize the volume fraction while considering only the fundamental buckling load constraint (i.e., no stress constraint). The force is uniformly distributed across four nodes at the midpoint of the upper edge, with a total force magnitude of $|F| = 10^{-3}$. The buckling constraint, $\mu^*$ *is set at* $\mu^*$ =0.15. The B-spline knot span length is selected as Δ=0.01L, where L represents the length of the square's sides. Fig. 4 illustrates the optimized design obtained using these parameters. The convergence history of the objective function (volume fraction) and the left-hand side of the buckling constraint is depicted in Fig. 5. The volume fraction decreases consistently until convergence is reached. It should be noted that the computation of the fundamental buckling load requires solving an eigenvalue problem, which is highly sensitive to boundary movement. As a result, some local small fluctuations can be observed during the optimization process. Nevertheless, the optimization leads to a significant reduction in the volume fraction, which decreases from an initial value of 0.61 to 0.201 after 500 iterations. This example demonstrates the effectiveness of the B-spline velocity field level set topology optimization approach in addressing the buckling constraint. By focusing solely on the fundamental buckling load constraint, the optimization process successfully minimizes the volume fraction while satisfying the constraint. The resulting design showcases the importance of selecting appropriate parameters and constraints to achieve an optimized structure that meets the desired performance criteria.

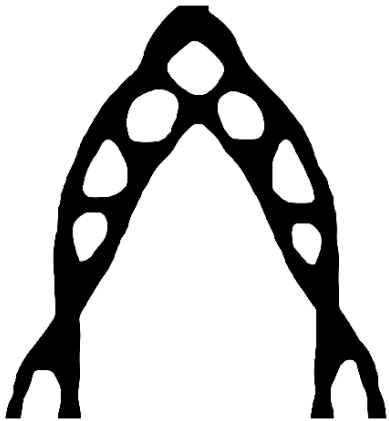

Figure 4. optimized design.

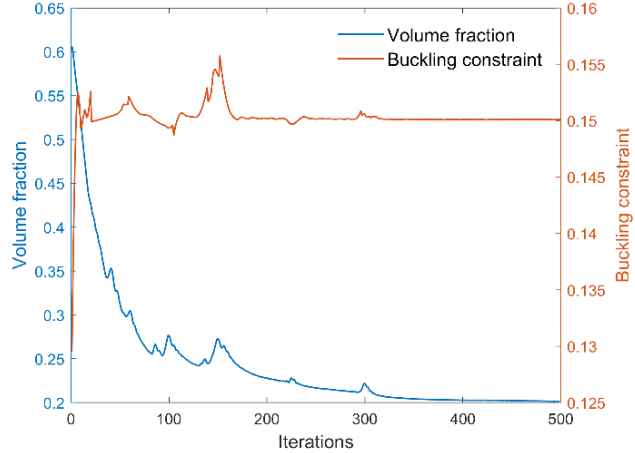

Figure 5. Convergence history of buckling constraint design.

**4.1.3 Volume minimization with buckling and stress constraint**

In this example, the optimization process incorporates both volume minimization with stress and buckling load constraints. Similar to the previous examples, the p-norm stress constraint, $\sigma^*$, and the buckling constraint, $\mu^*$, are set at $\sigma^* = 1.3$ and $\mu^* = 0.15$, respectively. The loading force on the mid-upper surface for stress and buckling constraints is consistent with those applied in sections 4.1.1 and 4.1.2. The optimized design can be observed in Fig. 6(a), while the von Mises stress distribution is illustrated in Fig. 6(b). The evolution of the material layout throughout the optimization process is depicted in Fig. 7. Notably, due to the presence of the buckling load constraint, the final design does not feature slender beam structures, as these would be more susceptible to buckling. Fig. 8 presents the convergence history, highlighting that the volume fraction decreases steadily from an initial value of approximately 0.6 to a final value of 0.216 during the optimization process. This reduction in volume fraction indicates that the optimization has been successful in minimizing material usage while ensuring that both stress and buckling constraints are satisfied. This example demonstrates the effectiveness of combining both stress and buckling load constraints in the optimization process. By addressing both constraints simultaneously, the resulting optimized design achieves a balance between stress resistance and buckling prevention, leading to a more effective structure.

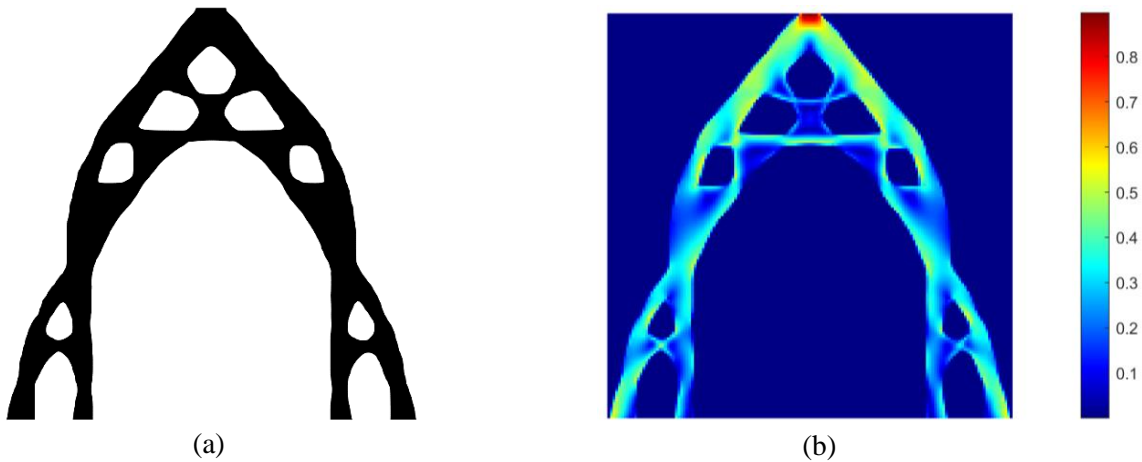

(a)  (b)

Figure 6. Optimized result: (a) material layout; (b) von Mises stress distribution.

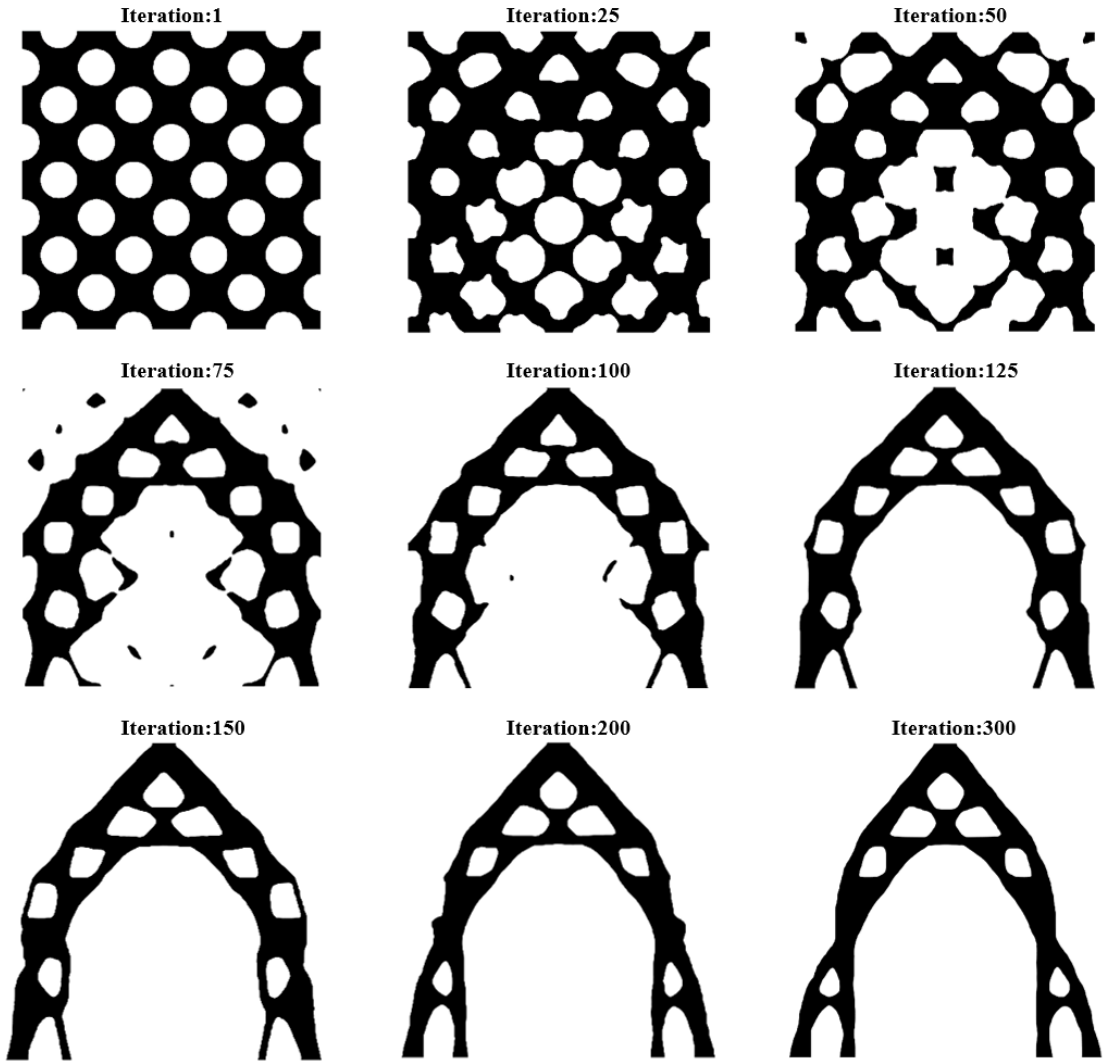

Figure 7. Evolution of material layout.

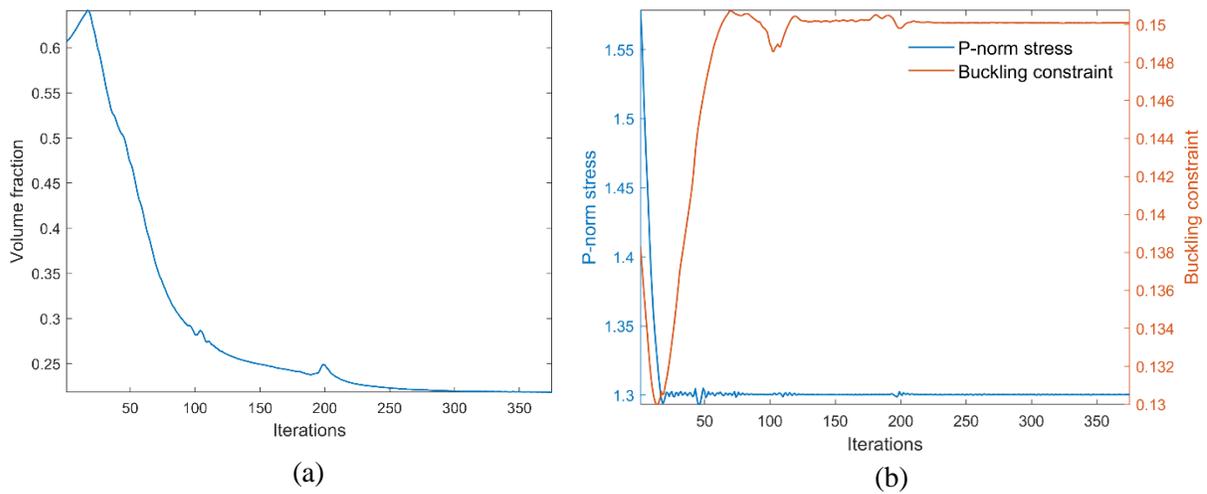

Figure 8. Convergence history. (a) objective function; (b) constraints.

## 4.2 L-bracket design

The dimensions of the L-bracket are illustrated in Fig. 11. It is important to note that the domain is discretized using a $100 \times 100$ finite element mesh, from which a $60 \times 60$ section is removed to create the L-bracket domain. A vertical load of $F_p = 1$ is applied to the upper right corner for the stress constraint problem. In this case, the elastic modulus and Poisson's ratio are chosen as $E = 1$ and $\mu = 0.3$, respectively. The B-spline velocity knots are uniformly distributed throughout the entire design domain, as depicted in Fig. 11. Assuming the maximum length of the L-bracket is represented by L, the value of the knot interval, $\Delta$, is selected as $\Delta = 0.02L$. Similarly, the bounds for the velocity design variables are set at ±0.2, and the moving limit of the MMA algorithm is chosen as 0.2.

### 4.2.1 volume minimization with stress constraint

In this example, the volume minimization of the L-bracket model, as presented in Fig. 9, is examined. The P-norm stress constraint is chosen as $\sigma^* = 0.65$, and the p-norm parameter is set to $p = 10$. The initial and optimized designs of the L-bracket are illustrated in Fig. 10(a) and (b), respectively, while the von Mises stress distribution is depicted in Fig. 10(c). Notably, the final solution features rounded corners, resulting in a much smoother shape than the initial design. The final optimized von Mises stress is uniformly distributed throughout the design space, and the optimized stress distribution closely approximates a fully stressed design. The evolution of the material layout is demonstrated in Fig. 11. As seen in the figure, the boundaries of the internal holes shift and merge with one another, producing non-trivial and optimized shapes. The corresponding evolution of the von Mises stress distribution is plotted in Fig. 12. The optimization process takes approximately 500 iterations to converge, as shown in Fig. 13. The volume fraction of the optimized shape is 0.196, indicating a successful reduction in material usage. The final optimized structure approaches a fully stressed state, suggesting that the optimization process has effectively balanced stress distribution and material reduction. This example highlights the capabilities of the B-spline velocity field level set topology optimization approach in handling complex design problems, such as the L-bracket example.

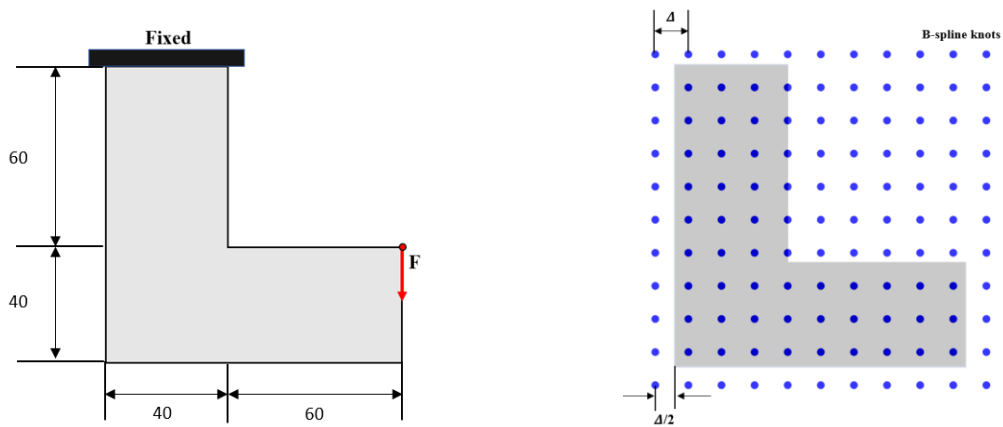

Figure 9. L-bracket design

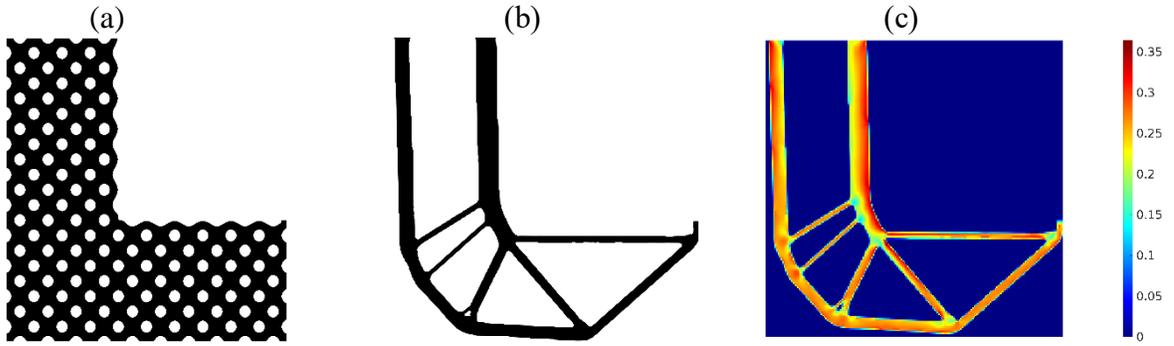

Figure 10. Volume minimization with p-norm stress constraint (a) initial design (b) optimized design (c) von Mises stress distribution

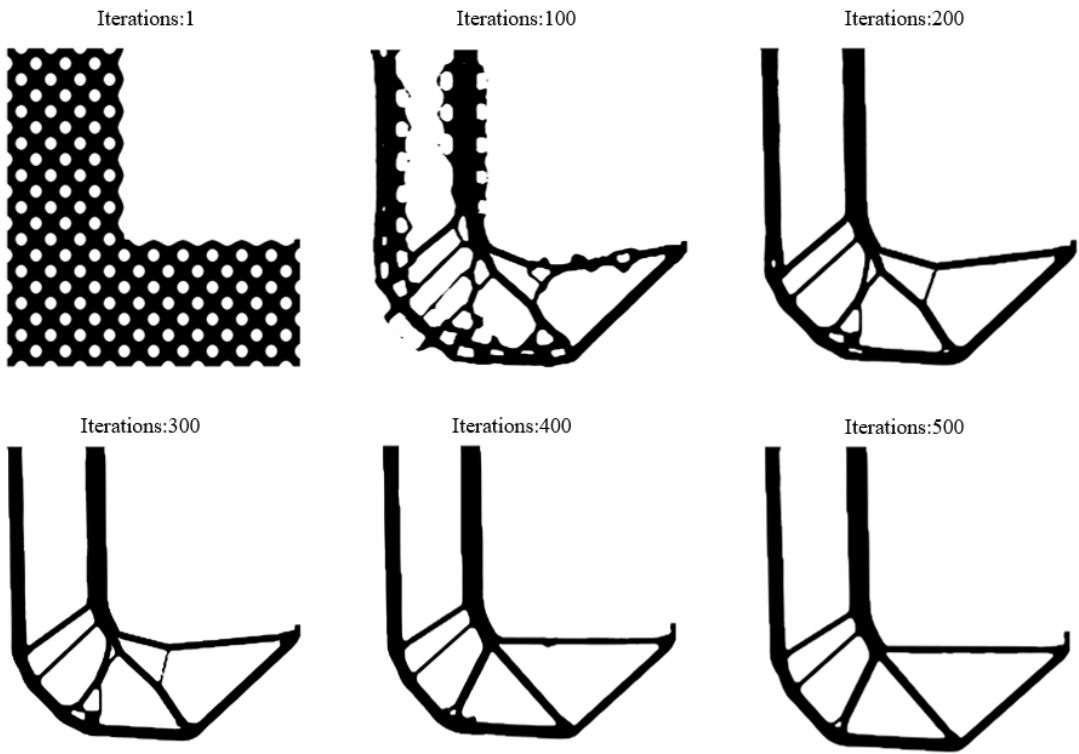

Figure 11. Evolution of material layout during optimization.

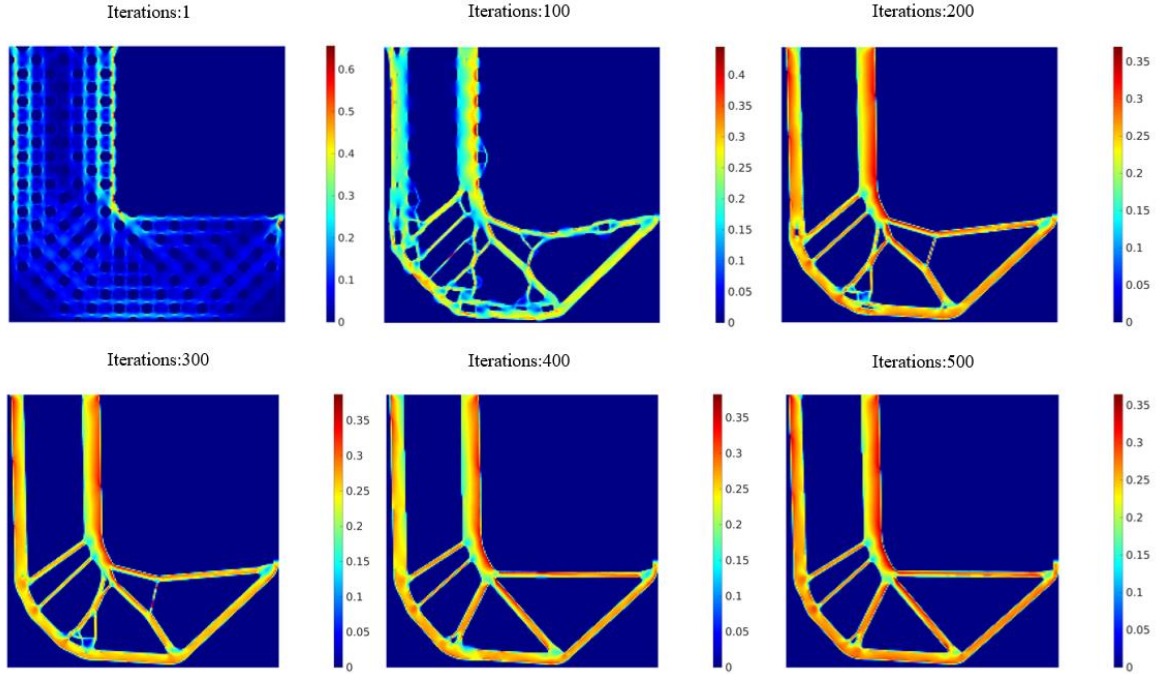

Figure 12. Evolution of von-Mises stress during optimization.

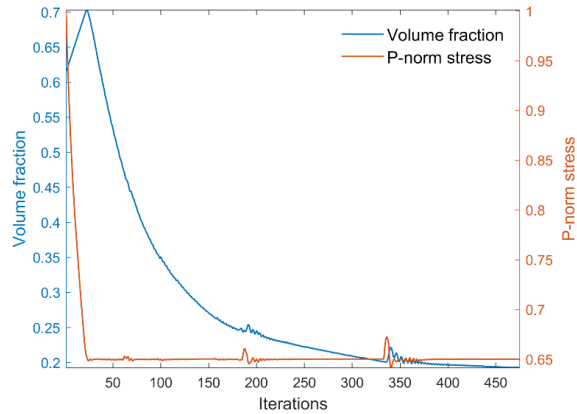

Figure 13. Convergence history.

### 4.2.2 Volume minimization with buckling constraint

In this section, the volume minimization of the L-bracket model with a buckling constraint is investigated using the proposed level set method. The finite element discretization and optimization configurations are consistent with those in section 4.2.1. A downward concentrated force of $F_b = 1 \times 10^{-3}$ is applied to the upper right corner. The buckling constraint, $\mu^*$, is selected as $\mu^* = 2.5$. Compared to the standard density-based method, local pseudo buckling modes in low-density regions [39] are not an issue for the level set method. The initial design and optimized design are illustrated in Fig. 14. In the final optimized design (Fig. 14(b)), no slender bars are present, while sharp corners are inevitably generated due to the absence of local stress constraints. It is evident that the optimized structural member in compression becomes wider to resist local buckling. The volume fraction reaches 0.307 after approximately 700 iterations. It is worth mentioning that local eigenvalues are highly sensitive to boundary movements and merges. Consequently, jumping phenomena of eigenvalues may occur during the optimization process. Therefore, local small fluctuations are observed in the convergence history, as shown in Fig. 15. This example demonstrates the capabilities of the level set method in handling optimization problems involving buckling constraints. The resulting optimized design features wider structural members in compression regions to resist local buckling, achieving a balance between volume minimization and buckling resistance. This highlights the importance of considering buckling constraints in the optimization process, particularly for structures subjected to compressive loads.

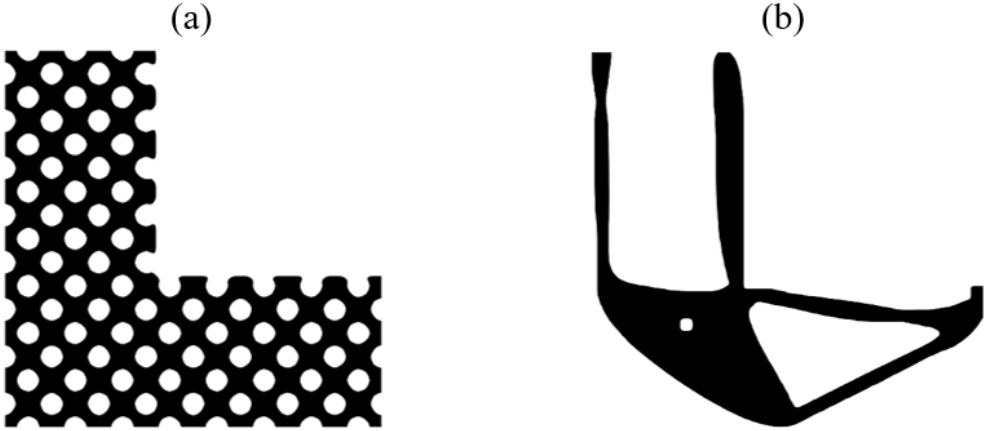

Figure 14. Volume minimization with buckling constraint (a) Initial design (b) Optimized design

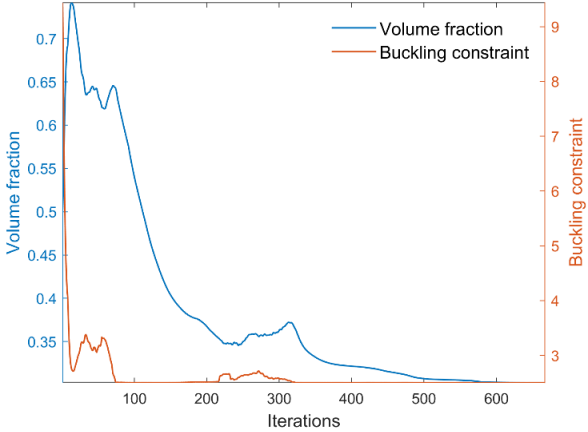

Figure 15. Convergence history

### 4.2.3 Volume minimization with stress and buckling constraint

In this section, the proposed methodology is applied to address the problem involving both stress and buckling constraints, focusing on the L-bracket design domain. The same material properties, discretization method, and optimization configuration from sections 4.2.2 and 4.2.1 are implemented. Two distinct loading cases are applied for stress and buckling analyses. For the stress constraint case, the load force is chosen as $F_p = 1$. For the buckling case, the loading is selected as $F_b = 1 \times 10^{-3}$. Similar to sections 4.2.2 and 4.2.1, the stress and buckling constraints are set as $\sigma^* = 0.65$ and $\mu^* = 2.5$.

Beginning with the initial design shown in Fig. 14(a), the solution for stress and buckling constraints is presented in Fig. 16(a), while the corresponding von Mises stress field is displayed in Fig. 16(b). As illustrated in Fig. 16(b), the maximum von Mises stress is evenly distributed near the rounded corner. When compared to the solution in section 4.2.2, the boundary in this case appears smoother, and no stress concentration points are found. The minimal volume fraction obtained is 0.509, indicating that the optimization has led to a reduction in material usage while satisfying both stress and buckling constraints. The optimization converges after approximately 500 iterations, and the evolution of topology shapes is demonstrated in Fig. 17. This example emphasizes the effectiveness of the proposed methodology in handling complex optimization problems involving multiple constraints. By considering both stress and buckling constraints, the resulting optimized design features smoother boundaries and more uniform stress distribution.

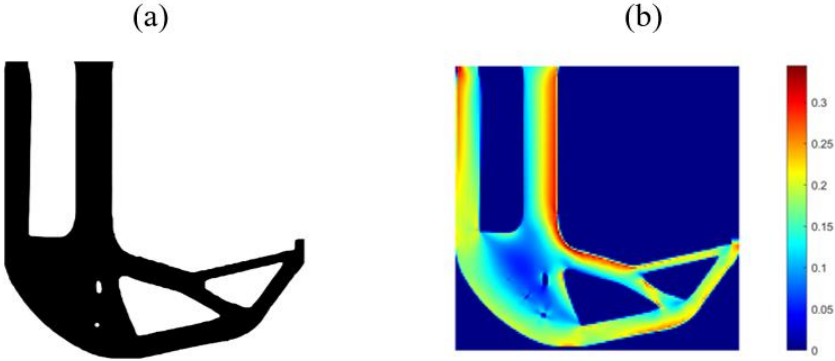

Figure 16. Volume minimization with stress and buckling constraint (a) Optimized design (b) von Mises stress distribution

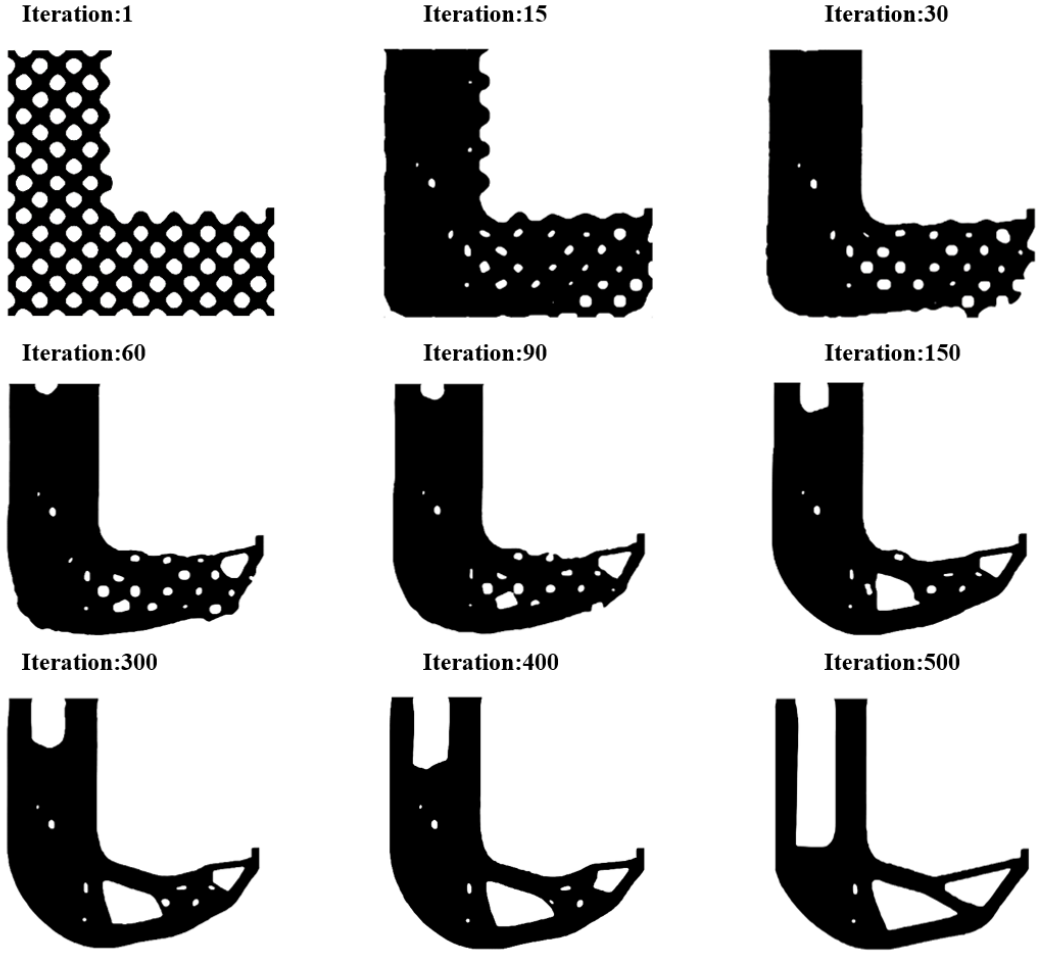

Figure 17. Evolution of material layout during optimization

## 4. Conclusion

In this paper, we introduce a novel velocity field level set (VFLS) method for topology optimization, which is based on the discrete adjoint method and employs B-spline representation to describe the velocity field. The analytical sensitivity of the proposed level set method can be fully computed using the discrete adjoint method. In this work, we derive and demonstrate the mathematical relationship between the level set velocity and the discrete adjoint sensitivities in detail. It is important to note that the presented sensitivity analysis method is not limited to B-spline space; any other velocity field representation can also be applied, as referenced in [40]. To showcase the effectiveness of our proposed sensitivity analysis method, we present stress and buckling constrained examples. These examples are typically considered challenging in the realm of level set topology optimization. The results demonstrate the effectiveness of the proposed method. Future work will focus on extending the current discrete adjoint-based level set method to address three-dimensional problems. Additionally, further research will investigate the integration of alternative velocity field representations and the development of more efficient algorithms to improve the computational efficiency of the method. This will ultimately expand the applicability of the discrete-adjoint based VFLS method to a broader range of engineering problems, making it a versatile tool for topology optimization.

Acknowledgments: The authors express their gratitude to Prof. Kazu Saitou for editing the manuscript.